\documentclass[11pt,oneside,reqno]{amsart}
\usepackage[utf8]{inputenc}
\usepackage[english]{babel}
\usepackage[babel]{csquotes}
\usepackage{graphicx}
\usepackage{ifthen}
\usepackage{stmaryrd}
\usepackage{amssymb}
\usepackage{amsmath}
\usepackage{amsthm}
\usepackage{fancyhdr}
\usepackage{array}
\usepackage{ragged2e}
\usepackage[colorlinks=true,linkcolor=black,citecolor=black,filecolor=black,urlcolor=black]{hyperref}

\allowdisplaybreaks

\def\R{\mathbb{R}}
\def\N{\mathbb{N}}

\def\T{\mathcal{T}}

\def\l{\mathit{l}}
\def\sf{f}
\def\sd{d}
\def\sx{x}

\def\stf{\mathcal{F}}
\def\std{d+1}
\def\stx{\mathcal{X}}
\def\stq{\mathcal{Q}}

\def\eps{\epsilon}

\def\ndofp{n_{\text{dof}}}

\newcommand{\tnorm}[1]{\left\vert\kern-0.9pt\left\vert\kern-0.9pt\left\vert #1\right\vert\kern-0.9pt\right\vert\kern-0.9pt\right\vert}

\DeclareMathOperator{\myspan}{span}

\DeclareMathOperator{\mydiv}{div}

\DeclareMathOperator{\const}{const}

\DeclareMathOperator{\sign}{sign}

\DeclareMathOperator{\conv}{conv}

\newcommand{\fIndex}[2][\empty]{%
  \ifthenelse{\equal{#1}{\empty}}
    {\emph{#2}\index{#2}}
    {\emph{#2}\index{#1}}}
\newcommand{\Index}[2][\empty]{%
  \ifthenelse{\equal{#1}{\empty}}
    {#2\index{#2}}
    {#2\index{#1}}}

\numberwithin{equation}{section}

\theoremstyle{plain} %
\newtheorem{theorem}{Theorem}[section]
\newtheorem{lemma}[theorem]{Lemma}
\newtheorem{corollary}[theorem]{Corollary}
\theoremstyle{definition} %
\newtheorem{definition}[theorem]{Definition}
\newtheorem{remark}[theorem]{Remark}
\theoremstyle{plain} %

\begin{document}
\title{
$L^\infty(L^\infty)$-boundedness of DG($p$)-solutions for nonlinear conservation laws
with boundary conditions
}

\author{Lutz Angermann}
\address{University of Technology at Clausthal, Department of Mathematics,\newline
Erzstrasse 1, D-38678 Clausthal-Zellerfeld, Germany}
\email{angermann@math.tu-clausthal.de}

\author{Christian Henke}
\address{University of Technology at Clausthal, Department of Mathematics,\newline
Erzstrasse 1, D-38678 Clausthal-Zellerfeld, Germany}
\email{henke@math.tu-clausthal.de}

\date{\today}

\pagestyle{fancy}
\rhead{$L^\infty(L^\infty)$-boundedness of DG($p$)-solutions \quad \thepage}
\lhead{}
\chead{}
\rfoot{}
\lfoot{}
\cfoot{}
\renewcommand{\headrulewidth}{0pt}

\begin{abstract}
We prove the $L^\infty(L^\infty)$-boundedness of a higher-order shock-cap\-turing
streamline-diffusion DG-method based on polynomials of degree $p\geq 0$
for general scalar conservation laws.
The estimate is given for the case of several space dimensions
and for conservation laws with initial and boundary conditions.
\end{abstract}

\maketitle
 
\section{Introduction}
In this paper we extend the analysis of a shock-capturing streamline-diffusion
DG-method for hyperbolic conservation laws in several space dimensions
which goes back to
\cite{Jaffre.Johnson.ea_Convergenceofdiscontinuous_1995}.
The original DG-method from \cite{Jaffre.Johnson.ea_Convergenceofdiscontinuous_1995}
is based on polynomials of maximal degree $p \geq 0$ (DG($p$)-method)
and is applied to a pure Cauchy problem.
Here we formulate the method for scalar conservation laws
with initial and boundary conditions.
To describe the further features of the method,
we recall the following sufficient conditions for convergence of a sequence
of approximate solutions \cite[Remark 1.2]{Szepessy_Measure-valuedsolutionsof_1989}: 
\begin{enumerate}
\item uniform boundedness in the $L^\infty(L^\infty)$-norm, i.e. $L^\infty$
in time and $L^\infty$ in space,
\item weak consistency with all entropy inequalities,
\item strong consistency with the initial condition.
\end{enumerate}
In the case of an unbounded domain, the condition $(1)$ can be replaced by the
\begin{enumerate}
\item[(1*)] uniform boundedness in the $L^\infty(L^2)$-norm,
\end{enumerate}
which was done in \cite{Jaffre.Johnson.ea_Convergenceofdiscontinuous_1995}. 
To the best of our knowledge, the only attempt to prove
the $L^\infty(L^\infty)$-boundedness
without using a finer auxiliary triangulation is given in
\cite{Johnson.Szepessy.ea_convergenceofshock-capturing_1990}
for the case $p=1.$
This proof can be extended for $p>1$ if the shock-capturing terms are defined
on finer triangulations \cite{Szepessy_Convergenceofstreamline_1991}.
Thanks to $(1^*)$ this is not necessary for the DG($p$)-method in
\cite{Jaffre.Johnson.ea_Convergenceofdiscontinuous_1995}. 
Our result presented here uses the skeleton of the proof from
\cite{Szepessy_Convergenceofstreamline_1991}, which is based on choosing
the test functions $v=I_h^p(U^{q-1})$ with a large even number $q$,
where $I_h^p$ is the Lagrange interpolation operator and
$U$ denotes the approximate solution.
Within this proof we use a new algebraic argument to verify the coercivity
of the shock-capturing term when $v=I_h^p(U^{q-1}).$

Let us recall the key points of the DG($p$)-method.
First, we have to choose the numerical flux on the element boundaries.
In contrast to \cite{Jaffre.Johnson.ea_Convergenceofdiscontinuous_1995},
where a strictly monotone numerical flux is necessary, we may also use
monotone numerical fluxes such as Engquist-Osher fluxes.
Second, there are two stabilization mechanisms. The DG($p$)-method
from the last mentioned reference contains a streamline-diffusion term
and a residual-based shock-capturing term.
The first term adds an anisotropic artificial viscosity and the second one
introduces some isotropic artificial viscosity locally
where the solution is nonsmooth.

The paper is organized as follows: In Section 2, we prepare some basic material
on scalar hyperbolic conservation laws with initial and boundary conditions.
Then, in Section 3 we introduce the DG($p$)-method under consideration,
and in Section 4 the condition (1*) is verified. 
After this we present our main Theorem \ref{Th:Linfty_main}
which is proved in Sections 5 and 6.
Here we give some background material on spatial and algebraic numerical ranges
and extend the condition (1) to the case $p>1.$

The entropy consistency and the consistency with the initial condition
will be proved in a forthcoming paper.

\section{Hyperbolic conservation laws with Boundary Conditions}
Let $Q_T=(0,T) \times \Omega \subset \R^{\std}, T>0, \label{zrz}$ $d\in\N,$
be an open time-space domain with boundary
$\Sigma_T= ~ (0,T) \times \Gamma,\, \Gamma =\partial \Omega \label{rzrz} $
and with outward unit normal $n.$ In this time-space domain a point with position
$\sx=(x_1,x_2,\ldots,x_{\sd})^T$ at time $t=x_0$ has the coordinates $\stx=(x_0,\sx)^T.$ 
Standard notation is used for the space of functions of bounded variations $BV(Q_T),$
Lebesgues spaces $L^q(Q_T)$ and Sobolev spaces $W^{l,q}(Q_T)$, $l \in \N$,
$1 \leq q \leq \infty.$

We consider for $u:Q_T \to \R$ the initial-boundary value problem
\begin{align}
L(u)=\nabla \cdot  \stf(u)  &=0 \text{ in } \, Q_T,
\label{hyper_problem} \\
u(0,\cdot)&=u_0 \text{ on } \, \Omega,
\label{hyper_problem_ic}
\end{align}
with the following boundary condition: For all $k \in \R, r \in \Sigma_T$
\begin{equation}
 \left(\sign(\gamma u(r) -k) 
- \sign(g_D(r)-k)\right) \left( \sf(\gamma u(r)) -\sf(k) \right) \cdot n(r) \geq 0,
\label{hyper_problem_bc}
\end{equation}
where  $\stf=(\cdot,\sf)^T:\R \to \R^{\std},$ $u_0: \Omega \to \R,$ $g_D:\Sigma_T \to \R$
are given smooth functions and $\gamma: Q_T \to \Sigma_T$ denotes a trace operator.
The function $\sign: \R \to \R$ is defined by
\begin{equation*}
\sign(x) = 
\begin{cases}
x/|x|, &\quad x \neq 0, \\
0,     &\quad x=0.
\end{cases}
\end{equation*}
Due to the hyperbolic nature of $(\ref{hyper_problem}),$ a boundary condition
of the form $u=g_D$ on $\Sigma_T$ usually over-determines the problem.
The generalization of the inflow boundary condition (where $\sf'(g_D) \cdot n \leq 0$)
for nonlinear $\sf$ also leads to a problem that is not well-posed.
This difficulty does not occur in $(\ref{hyper_problem})-(\ref{hyper_problem_bc}),$
because the solution $u_\eps$ of 
\begin{equation}
\begin{split}
  - \epsilon \Delta u_\eps + L(u_\eps) &=0 \text{ in } \, Q_T,\\
  u_\eps&=g_{\eps D} \text{ on } \, \Sigma_T,\\
  u_\eps(0,\cdot)&=u_{\eps 0} \text{ on} \, \Omega,
\label{para_problem}
\end{split}
\end{equation}
converges a.e.\ to a function $u \in BV(Q_T)$ as $\eps \to 0$ , which satisfies 
$(\ref{hyper_problem})-(\ref{hyper_problem_bc})$
\cite[Theorem 1]{Bardos.Roux.ea_Firstorderquasilinear_1979}.
Moreover, we can use this vanishing-viscosity method even in the space
$L^{\infty}(Q_T).$ It is possible to define a well-posed initial-boundary value problem,
which admits a unique weak entropy solution $u \in L^{\infty}(Q_T)$
\cite[Definition 7.2, Theorems 7.28, 8.20]{Malek.Nevcas.ea_WeakandMeasure-valued_1996}.

Let us shortly recall the concept of entropy pairs.
We say that $\stq=(\eta,q_1,\ldots,q_d) $ is an entropy pair
if $\eta : \R \to \R$ is continuous and convex, the entropy flux $q_j:\R \to \R$
is continuous and $\eta,q_1,\ldots,q_d$ satisfy for all $u \in \R$ the compatibility condition
\begin{equation}
\eta'(u) f_j'(u)= q_j'(u). 
\label{comp_bed}
\end{equation} 
For scalar conservation laws this is trivially fulfilled if the entropy flux
is defined as
\begin{equation}
q_j(u)=\int_{g_D}^u \eta'(r) f_j'(r)\, dr, \quad 1\leq j \leq \sd.
\label{def_q}
\end{equation}
 
\section{Formulation of the DG($p$)-method}
In this section we introduce the DG($p$)-method. To discretize
$(\ref{hyper_problem})-(\ref{hyper_problem_bc}),$ let 
$Q_{n,n+1}=(t_n,t_{n+1}) \times \Omega, \,Q_n= \{t_n\} \times \Omega \label{pzrz} \label{rpzrz}$
for the sequence of discrete time levels $0=  t_0<t_1<\dots<t_N, N \in \N,$
be a time-space decomposition of $Q_T.$
The boundary is defined by $\Sigma_{n,n+1}=(t_n,t_{n+1}) \times \Gamma \label{przrz}$
and $\Sigma_n=\{t_n\} \times \Gamma  \label{rprzrz}$. 

Consider an affine decomposition $\T_h^n$ of $Q_{n,n+1}$ belonging to a family
of quasi-uniform, admissible decompositions of $Q_{n,n+1},$ cf.   
\cite[Definition 1.49, 1.53, 1.140]{Ern.Guermond_Theoryandpractice_2004} 
into simplices or quadrilaterals $T$ and write $\T_h=\bigcup_{n \geq 0} \T_h^n.$
Let $h_T$ be the diameter of $T$ and $h$ the maximal diameter of all $T \in T_h^n.$

Set
\begin{equation}
W_h^n=\left\{ w \in L^2(Q_{n,n+1}): w|_T \in \mathbb{P}_p(T)
\quad \forall T \in \T_h^n \right\}, \quad W_h=\prod_{n\geq0} W_h^n,
\label{fem_space}
\end{equation}
where 
\begin{equation*}
\mathbb{P}_p(T)=\underset{\alpha \in \N_0^{\std}, \, |\alpha| \leq p}{\myspan}
\left\{ \stx^\alpha \right\},\, \stx \in T
\end{equation*}
is the space of polynomials of maximal degree $p$ defined on $T.$

We are now ready to define the DG($p$)-method and introduce by $R_{n,n+1}^i,$ $R_n^i$
the set of all interior faces of $Q_{n,n+1},$ $Q_n$ and by 
$\Lambda_{n,n+1},$ $\Lambda_n$ the set of all boundary faces.
We further set $R_{n,n+1}=R_{n,n+1}^i \cup \Lambda_{n,n+1}$ and
$R_{n}=R_{n}^i \cup \Lambda_{n}.$
In order to be able to describe discontinuous functions,
we denote by $\tau$ the common face shared by the elements $T^+=T$ and $T^-.$
We also define the normal vectors $n_T=n_T^+$ and $n_T^-$ on $\tau.$
Then we introduce the notation
\begin{equation}
v^\pm(x)=\underset{\mu \to +0}{\lim} v(x - \mu n^\pm),
\quad v_\pm^n(x)=v(t_n \pm 0,x_1,\ldots,x_{\sd})
\label{function_pm}
\end{equation}
and
\begin{equation}
\{v\}=\frac{1}{2} (v^+ + v^-), \quad \llbracket v n \rrbracket=v^+ n^+ + v^- n^- ,
\quad \llbracket v \rrbracket =v^+ - v^- .
\label{def_avg_jump}
\end{equation}
If $v$ is a function on $\Lambda_{n,n+1}$ or $\Lambda_{n},$ we write $v^-=g_D.$
By $C$ we will denote a positive constant independent of $h,$ not necessary
the same at each occurrence. 

Introducing the bivariate form
\begin{equation}
a(v,w)=\sum_{n=0}^{N-1}\sum_{T\in \T_h^n} \left\{ \int_T L(v) w \, d\stx
+ \int_{\partial T} (\hat{\stf}(v) -\stf(v^+) \cdot n^+) w^+ \, ds \right\},
\label{def_a}
\end{equation}
the DG($p$)-method for $(\ref{hyper_problem})-(\ref{hyper_problem_bc})$
can be formulated now: Find $U \in W_h$ such that for
$n=0,1,\cdots,N,$ $U \equiv U_{Q_{n,n+1}} \in W_h^n$ and for all $v \in W_h^n$
\begin{equation}
 a(U,v) +\sum_{T\in \T_h^n} \left\{ \left(\delta  L(U), \stf'(U) \cdot \nabla v \right)_{0,T}
+ \hat{\eps}  \left(\nabla U, \nabla v \right)_{0,T} \right\}=0,
\label{weak_problem_jaffre}
\end{equation}
where
\begin{align*}
\delta&=\delta(U)=C_1 h_T \left( \|\stf'(U) \|_{\l^2} \right)^{-1}, \\
\hat{\eps}&= \hat{\eps}(U)=\max\left( C_2 h^{2-\beta} R(U),C_3 h^{p+1/2} \right),
\ 0 < \beta < \frac{1}{2},\\
R(U)|_{T} &= \underset{T}{\max} \left(| L(U)|\right) +  
\frac{1}{h_T} \left( \underset{\partial^* T}{\max} \left(\left|\llbracket \stf(U) n\rrbracket \right|\right) 
+\underset{\partial^* T}{\max} \left(C_T \left| \llbracket U \rrbracket \right|\right) \right), \\
\partial^* T&=\{x \in \partial T: x \not\in R_{n+1}\}
\end{align*}
and $C_1,C_2,C_3>0.$ As mentioned in the introduction, $(\ref{weak_problem_jaffre})$
contains a streamline-diffusion term and a residual-based shock-capturing term.
Due to the $h$-dependency of this term, the DG-method can be interpreted as
a discrete vanishing-viscosity method.

The numerical flux $\hat{\stf}(U)$ is given by
\begin{equation}
\hat{\stf}(U)=\left\{ \stf(U)\right\} \cdot n^+ +C_T(U^+,U^-,n^+)  \llbracket U \rrbracket
\label{lf}
\end{equation}
and 
\begin{equation}
C_T(v^+,v^-,n^+)=
\begin{cases}
\frac{1}{2} &\quad n^+= \pm(1,0,\dots,0),\\
C_0^{\partial \Omega}\geq \int_0^1 |\stf'(v^- +s \llbracket v \rrbracket) \cdot n^+| \, ds &\quad v^-=g_D, \\
C_0^\Omega \geq \frac{1}{2}
\int_0^1 |\stf'(v^- +s \llbracket v \rrbracket ) \cdot n^+|\, ds &\quad \text{otherwise}.
\end{cases}
\label{ct_bed}
\end{equation}
Within this framework there are the following well-known numerical fluxes:\\
The Engquist-Osher flux if $f(0)=0$:
\begin{equation}
C_T(v^+,v^-,n^+)=
\begin{cases}
C_0^{\partial \Omega}=  \int_0^1 |\stf'(v^- +s \llbracket v \rrbracket ) \cdot n^+|\, ds    &\quad v^-=g_D, \\
C_0^\Omega = \frac{1}{2} C_0^{\partial \Omega}     &\quad \text{otherwise}.
\end{cases}
\label{def_EQ}
\end{equation}
The Lax-Friedrichs flux:
\begin{equation}
C_T(v^+,v^-,n^+)=
\begin{cases}
C_0^{\partial \Omega}=  \underset{z\in[v^+,v^-]}{\sup}|\stf'(z) \cdot n^+|     &\quad v^-=g_D, \\
C_0^\Omega = \frac{1}{2} C_0^{\partial \Omega}     &\quad \text{otherwise}.
\end{cases}
\label{def_LF}
\end{equation}
Further $C_0$ is a positive constant satisfying
\begin{equation}
C_0 \geq \|\stf'\|_{0,\infty,\R}= \underset{x \in \R}{\max} \|\stf'(x)\|_{\l^2} 
\label{r01}
\end{equation}
and
\begin{equation}
C_0^{\partial \Omega},\, C_0^{\Omega} \leq C_0.
\label{r02}
\end{equation}

\section{Preliminaries}
In this section we want to verify the uniform $L^\infty(L^2)$-bound (1*)
which is based on \cite{Jaffre.Johnson.ea_Convergenceofdiscontinuous_1995}.
As a byproduct of the study of this technical result, we introduce some notation
and prepare the basic material for the $L^\infty(L^\infty)$-bound
presented in Section \ref{sec:infinfbnd}.
To make this precise, we choose $\eta(U)=U^2/2,$ $\varphi=1$ and $v=\eta'(U) \varphi.$
The main theorem of this paper can be obtained by applying the test function
$v=I_h^p(\eta'(U) \varphi),$ where $\eta(U)= U^q/q,$ $\varphi=1$ and
$I_h^p:C(\overline{Q_T}) \to W_h$ is the Lagrange interpolation operator.
Here, due to the fact that $\eta'(U) \notin W_h,$ an interpolation or projection operator
is necessary. Notice that in this case we get an additional difficulty to estimate
terms which contain the difference $\eta'(U) \varphi-I_h^p(\eta'(U) \varphi).$

By the definition of the bivariate form
\begin{equation}
b(v,w)=a(v,w)+\sum_{n=0}^{N-1} \sum_{T \in \T_h^n}  
\left( \delta(v) L(v), \stf'(v) \cdot \nabla w \right)_{0,T} 
\label{b1}
\end{equation}
we have, for $w=\eta'(U) \varphi$ and an entropy pair $(\eta,q)$
satisfying $(\ref{comp_bed})$ and $(\ref{def_q}),$ that
\begin{equation}
\begin{split}
&b(v,\eta'(v) \varphi)=\sum_{n=0}^{N-1} \sum_{T \in \T_h^n} \bigg\{ 
 \left( \delta(v) L(v), \stf'(v) \cdot \nabla (\eta'(v) \varphi) \right)_{0,T} \\
&+  \int_{\partial T} \frac{1}{2} \llbracket \stf(v) n \rrbracket \eta'(v)\varphi \, ds
+ \int_{\partial T} C_T \llbracket v \rrbracket \eta'(v) \varphi \, ds \bigg\} \\
&+\int_{\partial T}  \stq(v) \cdot n \varphi \, ds 
- \int_{[0,t_N]\times \Omega}  \stq(v) \cdot \nabla \varphi \, d\stx.
\end{split}
\label{b_3}
\end{equation}
As usual in DG-methods we consider the different behaviour of inner and boundary faces
\begin{equation}
\begin{split}
b(v,\eta'(v) \varphi) &=
\int_\Omega \eta(v_-^N) \varphi^N \, d\sx - \int_\Omega \eta(v_-^0) \varphi^0 \, d\sx \\
&- \int_{[0,t_N]\times \Omega}  \stq(v) \cdot \nabla \varphi \, d\stx\\
&+\sum_{i=0}^{5}E_i(f,\eta,v,\varphi)- F(f,\eta,v,\varphi),
\end{split}
\label{b1_2}
\end{equation}
where
\begin{align}
E_0(f,\eta,v,\varphi)&=\sum_{n=0}^{N-1} \sum_{T \in \T_h^n} 
 \left( \delta(v)L(v), \stf'(v) \cdot \nabla (\eta'(v) \varphi)\right)_{0,T},
\label{e0} \\
E_1(f,\eta,v,\varphi)&=\sum_{n=0}^{N-1} \int_\Omega \left( \eta(v_-^n)-\eta(v_+^n)
- \eta'(v_+^n)(v_-^n -v_+^n) \right) \varphi^n \, d\sx,
\label{e1} \\
E_2(f,\eta,v,\varphi)
& =\sum_{n=0}^{N-1} \sum_{\tau \in R_{n,n+1}^i} \int_\tau \left( \llbracket  \stq(v) n\rrbracket
-  \llbracket \stf(v) n\rrbracket \left\{ \eta'(v)\right\} \right) \varphi \, ds,
\label{e2} \\
E_3(f,\eta,v,\varphi)
&=\sum_{n=0}^{N-1} \sum_{\tau \in R_{n,n+1}^i}
\int_\tau C_0^{\Omega} \int_0^1 \eta''(v^- +r\llbracket v \rrbracket)dr\llbracket v \rrbracket^2   \varphi \, ds,
\label{e3} \\
E_4(f,\eta,v,\varphi)
& =\sum_{n=0}^{N-1} \sum_{\tau \in \Lambda_{n,n+1}} \int_\tau \left( \llbracket  \stq(v)n \rrbracket
-  \llbracket \stf(v)n \rrbracket \left\{ \eta'(v)\right\} \right) \varphi \, ds,
\label{e4} \\
E_5(f,\eta,v,\varphi)
&=\sum_{n=0}^{N-1} \sum_{\tau \in \Lambda_{n,n+1}} \int_\tau  C_0^{\partial \Omega}
\int_0^1 \eta''(g_D +r\llbracket v \rrbracket)dr\llbracket v \rrbracket^2   \varphi \, ds, 
\label{e5}\\
F(f,\eta,v,\varphi)
&= -\sum_{n=0}^{N-1} \sum_{\tau \in \Lambda_{n,n+1}} \int_\tau \left( \frac{1}{2}\llbracket \stf(v) n \rrbracket
+C_0^{\partial \Omega} \llbracket v \rrbracket \right)\eta'(g_D)  \varphi \, ds.
\label{f}
\end{align}
Next, we will show the nonnegativity of $\sum_1^5 E_i(f,\eta,v,\varphi).$
By the convexity of $\eta,$ this is true for $E_1.$
In order to treat $E_2,$ we consider the expression
\begin{eqnarray*}
& &\left(  \stq(v^+) -  \stq(v^-) - \frac{1}{2} \left( \stf(v^+)  -\stf(v^-) \right) \left( \eta'(v^+)
+ \eta'(v^-) \right) \right) \cdot n^+ \\ 
&=& \left( \int_{v^-}^{v^+} \left( \stq' -\stf'\eta' \right) dr
+ \int_{v^-}^{v^+} \stf' \left( \eta' -\frac{1}{2} \left( \eta'(v^+)+\eta'(v^-) \right) \right) dr \right) \cdot n^+\\
&\underset{(\ref{comp_bed})}{=}& \int_0^1 \stf'(v^- + s\llbracket v \rrbracket)
\cdot n^+ \left( \eta'(v^- +s\llbracket v \rrbracket)
- \frac{1}{2}\left( \eta'(v^+) + \eta'(v^-) \right) \right) \llbracket v \rrbracket \, ds.\\ 
\end{eqnarray*}
By the properties of convex functions, it follows that
\begin{equation*}
\left(\eta'(v^- +s \llbracket v \rrbracket) - \eta'(v^-)\right) \llbracket v \rrbracket \geq 0
\end{equation*}
and
\begin{equation*}
\left(\eta'(v^- + s \llbracket v \rrbracket) - \eta'(v^+) \right) \llbracket v \rrbracket
= \left(\eta'(v^+ - (1-s)\llbracket v \rrbracket) - \eta'(v^+)\right) \llbracket v \rrbracket \leq 0.
\end{equation*}
Then we have
\begin{eqnarray*}
& & \left|\left( \eta'(v^- +s\llbracket v \rrbracket) - \frac{1}{2}\left( \eta'(v^+) + \eta'(v^-) \right) \right) \llbracket v \rrbracket \right|\\
&=& \left| \left( \frac{1}{2} \left( \eta'(v^- +s\llbracket v \rrbracket)
- \eta'(v^-)\right)+ \frac{1}{2}\left(\eta'(v^- +s \llbracket v \rrbracket)
- \eta'(v^+)  \right) \right) \llbracket v \rrbracket \right|\\
&\leq&  \frac{1}{2}\left| \left( \eta'(v^- +s\llbracket v \rrbracket)
- \eta'(v^-) \right)\llbracket v \rrbracket \right|
+ \frac{1}{2}\left|\left(\eta'(v^- +s \llbracket v \rrbracket)
- \eta'(v^+)  \right) \llbracket v \rrbracket \right|\\
&=& \frac{1}{2}\left( \left( \eta'(v^- +s\llbracket v \rrbracket)
- \eta'(v^-) \right)- \left( \eta'(v^- + s\llbracket v \rrbracket)
- \eta'(v^+)\right) \right)\llbracket v \rrbracket \\
&=& \frac{1}{2} \left( \eta'(v^+)-\eta'(v^-) \right) \llbracket v \rrbracket
= \frac{1}{2}  \int_0^1 \eta''(v^- +r\llbracket v \rrbracket  ) \, dr \llbracket v \rrbracket^2,\\
\end{eqnarray*}
which immediately implies that
\begin{eqnarray*}
& &\left |\left(  \stq(v^+) -  \stq(v^-)
- \frac{1}{2} \left( \stf(v^+)  -\stf(v^-) \right) \left( \eta'(v^+) + \eta'(v^-) \right) \right) \cdot n^+ \right| \\
&\leq& \int_0^1 \left|\stf'(v^- + s\llbracket v \rrbracket)
\cdot n^+\right| \left| \left( \eta'(v^- +s\llbracket v \rrbracket)
- \frac{1}{2}\left( \eta'(v^+) + \eta'(v^-) \right) \right) \llbracket v \rrbracket \right| \, ds\\ 
&\leq& \frac{1}{2}\int_0^1 \left|\stf'(v^- + s\llbracket v \rrbracket)
\cdot n^+\right| \int_0^1 \eta''(v^- +r\llbracket v \rrbracket  ) \, dr \llbracket v \rrbracket^2 \, ds\\ 
&\leq& \frac{1}{2}\int_0^1 \left|\stf'(v^- + s\llbracket v \rrbracket)
\cdot n^+\right|\,ds \int_0^1 \eta''(v^- +r\llbracket v \rrbracket  ) \, dr \llbracket v \rrbracket^2.\\ 
\end{eqnarray*}
Finally, having in mind the fact that $\eta''\geq 0,$ we can use 
the definitions $(\ref{e3})$ and $(\ref{ct_bed})$
\begin{equation}
\thinmuskip=0mu
\medmuskip=0mu
\thickmuskip=1mu
\begin{split}
&E_2(f,\eta,v,\varphi)+ E_3(f,\eta,v,\varphi)\\
&\geq \sum_{\overset{n=0}{\tau \in R_{n,n+1}^i }}^{N-1} \int_\tau \left( C_0^{\Omega}
-\frac{1}{2}\int_0^1 \left|\stf'(v^- + \xi \llbracket v \rrbracket) \cdot n^+\right|\,d\xi  \right) 
\left( \int_0^1  \eta''(v^- +r\llbracket v \rrbracket) \,dr \right) \llbracket v \rrbracket^2 \varphi \, ds\\
&\geq \;\;\;\; 0,
\end{split}
\label{e2e3}
\end{equation}
where we have used nonnegative test functions $\varphi.$
The same arguments as before lead to $E_4(f,\eta,v,\varphi)+ E_5(f,\eta,v,\varphi)\geq 0.$ 
\begin{remark}
The local condition $(\ref{ct_bed})$
\begin{equation*}
C_0^{\Omega} - \frac{1}{2}\int_0^1 \left|\stf'(v^- + \xi \llbracket v \rrbracket) \cdot n^+\right|\,d\xi \geq 0
\end{equation*}
for $(\ref{e2e3})$ allows a smaller constant $C_0^{\Omega}$
than \cite[Remark 2.5]{Jaffre.Johnson.ea_Convergenceofdiscontinuous_1995}
\begin{equation*}
C_0^{\Omega} - \frac{1}{2} \|\stf'\|_{0,\infty,\R}
\frac{\int_0^1 \int_0^s \left( n''(v^- +r \llbracket v \rrbracket)
+ \eta''(v^+ -r\llbracket v \rrbracket ) \right) \, dr \, ds}
{\int_0^1 \eta''(v^- +r \llbracket v \rrbracket) \, dr} \geq 0.
\end{equation*}
Thus, the corresponding numerical flux $\hat{\stf}(U)$ is a monotone flux function
($U^+ \mapsto \hat{\stf}(U)$ is increasing and $U^- \mapsto \hat{\stf}(U)$ is decreasing).
As mentioned in \cite{Jaffre.Johnson.ea_Convergenceofdiscontinuous_1995},
the requirement that $\hat{\stf}(U)$ is a strictly monotone numerical flux, e.g.
\begin{equation*}
C_0^{\Omega} - \frac{1}{2}\int_0^1 \left|\stf'(v^- + \xi \llbracket v \rrbracket)
\cdot n^+\right|\,d\xi \geq \epsilon >0,
\end{equation*}
is necessary for (2). More precisely, the condition
\begin{equation*}
\sum_{n=0}^{N-1} \sum_{T \in \T_h^n} \int_T h R(U)^2 \, d\stx \leq C \|u_0\|_{0,2,\Omega}^2
\end{equation*}
has to be fulfilled.
\end{remark}
Let us now consider the equation $(\ref{weak_problem_jaffre})$
with $v=\eta'(U) \varphi$:
\begin{equation}
\begin{split}
&\quad \; b(U,\eta'(U) \varphi) + \sum_{n=0}^{N-1} \sum_{T \in \T_h^n} 
\hat{\eps}(U) \left(\nabla U,\nabla( \eta'(U) \varphi) \right)_{0,T} \\
&= \int_\Omega \eta(U_-^N) \varphi^N \, d\sx - \int_\Omega \eta(U_-^0) \varphi^0 \, d\sx 
+ \sum_{n=0}^{N-1} \sum_{T \in \T_h^n} \hat{\eps}(U) \left(\nabla U,\nabla( \eta'(U) \varphi) \right)_{0,T}\\
&- \int_{[0,t_N]\times \Omega}  \stq(U) \cdot \nabla \varphi \, d\stx
+\sum_{i=0}^{5}E_i(f,\eta,U,\varphi)- F(f,\eta,U,\varphi)=0.
\end{split}
\label{weak_problem2}
\end{equation}
Therefore, in the case $\eta(U)=U^2/2$ and $\varphi=1,$ we obtain
\begin{equation}
\begin{split}
&\frac{1}{2} \int_\Omega (U_-^N)^2 \, d\sx 
 + \sum_{n=0}^{N-1} \sum_{T \in \T_h^n} 
\hat{\eps}(U) \left(\nabla U,\nabla U \right)_{0,T} \\
&+\sum_{i=0}^{5}E_i(f,U^2/2,U,1) = \frac{1}{2} \int_\Omega (u_0)^2  \, d\sx +F(f,U^2/2,U,1),
\end{split}
\label{weak_problem_uu}
\end{equation}
where we have used $U_{-}^0=u_0.$ By the help of Young's inequality we see that
\begin{equation}
\begin{split}
F(f,U^2/2,U,1) 
&\leq \sum_{n=0}^{N-1} \sum_{\tau \in \Lambda_{n,n+1}}
\int_\tau \frac{1}{2} C_0^{\partial \Omega}  \llbracket U \rrbracket^2 
+ \frac{9}{8} C_0^{\partial \Omega} g_D^2  \, ds \\
&= F_1(f,U^2/2,U,1)+F_2(f,g_D^2/2,g_D,1),
\end{split}
\label{rhs_L2}
\end{equation}
thus we arrive at $E_4(f,U^2/2,U,1)+E_5(f,U^2/2,U,1)-F_1(f,U^2/2,U,1) \geq 0$
and
\begin{equation}
 \frac{1}{2} \|U_-^N\|^2_{0,2,\Omega}
+ \sum_{n=0}^{N-1} \sum_{T \in \T_h^n}  \|\delta(U)^{1/2} L(U)\|^2_{0,2,T} 
\leq \frac{1}{2} \|u_0\|^2_{0,2,\Omega} 
+ \frac{9}{8}  C_0 \|g_D \|_{0,2,\Sigma_T}^2. 
\label{weak_problem_uu4}
\end{equation}
Finally, using  for $t_{N-1} \leq t \leq t_N$ and $q \in \N$ the identity
\begin{equation}
 \|U(t, \cdot)\|_{0,q,\Omega}^q 
=  \|U_-^N\|_{0,q,\Omega}^q - q \int_t^{t_N} \int_\Omega U^{q-1}(t',\sx) \, \mydiv \, \stf(U(t',\sx)) \, d\sx dt',
\label{id_Uq}
\end{equation} 
Young's inequality yields
\begin{equation}
\begin{split}
 \|U(t, \cdot)\|_{0,2,\Omega}^2
 &\leq \|U_-^N\|_{0,2,\Omega}^2
+ \frac{C_0}{2 C_1 \underset{T \in \T_h^{N-1}}{\min} \{h_T\} } \int_t^{t_N} \|U(t', \cdot ) \|_{0,2,\Omega}^2 \, dt' \\
&+ 2 \sum_{T \in \T_h^{N-1}}  \| \delta(U)^{1/2} L(U) \|_{0,2,T}^2   
\end{split}
\label{L_2abschaetzung3}
\end{equation}
and a Gronwall argument estimates the right-hand side of $(\ref{L_2abschaetzung3})$
by means of the left-hand side of $(\ref{weak_problem_uu4}).$
The quasi-uniformity of  $\left\{ \mathcal{T}_h^{N-1} \right\}_{h > 0}$
ensures the boundedness of the Gronwall constant and the following theorem results.
\begin{theorem}
Let $\Omega$ be a domain with a Lipschitz boundary and
$\left\{ \mathcal{T}_h \right\}_{h > 0}$ be a quasi-uniform family
of decompositions of $(0,T) \times \Omega.$
Let $U$ be a solution of $(\ref{weak_problem_jaffre})$
satisfying the assumptions $(\ref{r01})$ and $(\ref{r02}).$
Then there exists  a constant $C>0$ independent of $h,$ such that,
for all $t \in (t_{N-1},t_N),$
\begin{equation}
\begin{split}
\|U(t, \cdot)\|_{0,2,\Omega}
&\leq C \left( 
 \|u_0\|_{0,2,\Omega}
+  \|g_D \|_{0,2,\Sigma_T} 
\right).
\label{L_2abschaetzung4}
\end{split}
\end{equation} 
\label{Th:LinftyL2}
\end{theorem}

\section{$L^\infty(L^\infty)$-boundedness of DG($p$)-solutions}
\label{sec:infinfbnd}
As announced above, in this section we prove that
$\|U\|_{0,\infty,Q_T}$ is uniformly bounded.
The main idea is to control the interpolation error $U^{q-1}-I_h^p(U^{q-1})$
in the second argument of the bivariate form by the aid
of the special shock-capturing term.
At the end of this section we formulate two corollaries which are consequences
of the limiting process $h \to 0$ and of the special case $p=0,$ respectively.
Our main result is the following.
\begin{theorem}
Let $\Omega$ be a domain with a Lipschitz boundary and
$\left\{ \mathcal{T}_h \right\}_{h > 0}$ be a quasi-uniform family of decompositions
of $(0,T) \times \Omega.$ 
Let $U$ be a solution of $(\ref{weak_problem_jaffre})$
satisfying the assumptions $(\ref{r01})$ and $(\ref{r02}).$
Then there exists  a constant $C>0$ independent of $h$ such that
\begin{equation}
\|U \|_{0,\infty,Q_T} \leq C 
\left(\|u_0\|_{0,\infty,\Omega} + \|g_D \|_{0,\infty,\Sigma_T} + 1\right).
\end{equation}
\label{Th:Linfty_main}
\end{theorem}
The proof is based on the next lemma which contains the extension of
\cite[Lemma 3.3]{Szepessy_Convergenceofstreamline_1991} and
\cite[Lemma 4.2]{Szepessy_Convergenceofshock-capturing_1989}
for $p>1$ and which is proved in the last section.

\begin{lemma}
For Lagrange finite elements with a shape regular family of meshes
$\left\{ \mathcal{T}_h^n \right\}_{h>0}$
there is a constant $C >0$ independent of $q$ and $h$ such that
for all $v \in W_h$ and $q=2m,\, m \in \N:$
\begin{equation}
\left( \nabla v, \nabla I_h^p(v^{q-1})\right)_{0,T}
\geq C  \int_T \|\nabla v\|_{\l^2}^2 \|v\|_{0,\infty,T}^{q-2}\, d\stx,
\quad \forall T \in \mathcal{T}_h^n. 
\label{nabla_qm1_ungl}
\end{equation}
\label{nabla_qm1}
\end{lemma}
\begin{proof}[Proof of Theorem \ref{Th:Linfty_main}]
Setting $v=I_h^p(\eta'(U) \varphi)$ in $(\ref{weak_problem_jaffre})$
with $\eta(v)= v^q/q,$ $\varphi=1$ and $q>2$ an even natural number,
we obtain
\begin{equation}
\begin{split}
&\phantom{+\,} \frac{1}{q} \int_\Omega \left( U_-^N \right)^q \, d\sx 
+\sum_{i=0}^5 E_i(f,U^q/q,U,1) 
-\left( b(U,U^{q-1})-b(U,I_h^p(U^{q-1})) \right)\\
&+\sum_{n=0}^{N-1} \sum_{T \in \T_h^n} \hat{\eps}(U)
\left(\nabla U,\nabla I_h^p(U^{q-1})\right)_{0,T}
=\frac{1}{q} \int_\Omega \left( u_0 \right)^q d\sx
+ F(f,U^q/q,U,1).
\end{split}
\label{weak_problem4}
\end{equation}
Here the key point is that the interpolation error is bounded by
the isotropic shock-capturing term. To see this, a careful consideration
of the interpolation operator $I_h^p$ is necessary.

By 
\begin{equation*}
\int_{R_{n+1}}  \frac{1}{2} \llbracket \stf(U) n \rrbracket
+ C_T \llbracket U \rrbracket \, d\sx =0,
\end{equation*}
we express the interpolation error as
\begin{equation}
\begin{split}
&\phantom{=} \left( b(U,U^{q-1})-b(U,I_h^p(U^{q-1})) \right) \\
&\underset{(\ref{b1})}{=} \sum_{n=0}^{N-1} \sum_{T \in \T_h^n}
\Big\{ \int_T \nabla \cdot \stf(U) \left( U^{q-1}-I_h^p (U^{q-1}) \right) \, d\stx\\
&+  \int_T \delta(U) L(U) \stf'(U) 
\cdot \left( \nabla U^{q-1} - \nabla I_h^p\left( U^{q-1} \right) \right) \, d\stx \\
&+ \int_{\partial^* T} \frac{1}{2} \llbracket \stf(U)n \rrbracket
\left( U^{q-1} -I_h^p(U^{q-1}) \right) \, ds\\ 
&+ \int_{\partial^* T} C_T \llbracket U \rrbracket
\left( U^{q-1}-I_h^p(U^{q-1}) \right) \, ds \Big\}
=\sum_{n=0}^{N-1} \sum_{T \in \T_h} \sum_{i=1}^{4} A_T^i.
\end{split}
\label{int_err_q}
\end{equation}
Since $U|_T \in \mathbb{P}_p(T),$ we deduce that $|U|_{p+1,\infty,T}=0.$
Arguing as in \cite{Szepessy_Convergenceofstreamline_1991}[p. 765],
we may write for $q\geq 3$ 
\begin{equation}
|U^{q-1}|_{p+1,\infty,T}
\leq C q^{p+1} h^{-p+1} \| \nabla U \|_{0,\infty,T}^{2} \|U\|_{0,\infty,T}^{q-3}.
\label{diskret_hoch_q_minus_1_uneq}
\end{equation}
Together with a standard interpolation error, we thus conclude from the first term of
$(\ref{int_err_q})$ that
\begin{eqnarray*}
\left| A_T^1 \right| &\leq& \|(I-I_h^p) U^{q-1} \|_{0,\infty,T} \int_T |\nabla \cdot \stf(U) | \, d\stx \\
&\leq& C h_T^{p+1} |U^{q-1}|_{p+1,\infty,T} \int_T | \nabla \cdot \stf(U) | \, d\stx \\
&\underset{(\ref{diskret_hoch_q_minus_1_uneq})}{\leq}
& C q^{p+1} h_T^2 \|\nabla U\|_{0,\infty,T}^2 \|U\|_{0,\infty,T}^{q-3} \int_T |\nabla \cdot \stf(U) | \, d\stx\\
&\leq& C q^{p+1} h_T^2   \int_{T \cap\{|U| > 1  \}} \|\nabla U\|_{0,\infty,T}^2
\|U\|_{0,\infty,T}^{q-2} |\nabla \cdot \stf(U) | \, d\stx \\
&\quad+& C q^{p+1} h_T^2  \int_{T \cap\{|U| \leq 1  \}} \|\nabla U\|_{0,\infty,T}^2  |\nabla \cdot \stf(U) | \, d\stx  \\
&\leq& C q^{p+1} h_T^2  \max_T(|\nabla \cdot \stf(U) |) \|U\|_{0,\infty,T}^{q-2} \int_{T} \|\nabla U\|_{0,\infty,T}^2   \, d\stx \\
&\quad+& C q^{p+1} h_T^2  \max_T(|\nabla \cdot \stf(U) |) \int_{T} \|\nabla U\|_{0,\infty,T}^2  \, d\stx.\\
\end{eqnarray*}
By the quasi-uniformity of $\left\{ \T_h^n \right\}_{h>0}$
and an inverse inequality, we obtain that
\begin{equation}
\int_T \|\nabla v\|_{0,\infty,T}^2 \, d\stx \leq 
C \int_T \|\nabla v\|_{\l^2}^2 \, d\stx.
\label{int_norm_ungl}
\end{equation}
Thus we have, by Lemma \ref{nabla_qm1_ungl},
\begin{equation*}
|A_T^1|
\leq Cq^{p+1} h_T^2 \max_T(|\nabla \cdot \stf(U)|) \left\{ \left(\nabla U, \nabla I_h^p(U^{q-1})\right)_{0,T}
 + \|\nabla U\|_{0,2,T}^2 \right\}.
\end{equation*}
In a similar fashion, we can estimate the complete right-hand side of $(\ref{int_err_q}).$
Consequently, by $(\ref{weak_problem_uu4})$ we conclude that
\begin{eqnarray}
& &\left| \left( b(U,U^{q-1})-b(U,I_h^p(U^{q-1})) \right) \right| \notag \\
& \leq& C q^{p+1} \sum_{n=0}^{N-1} \sum_{T \in \T_h}  h_T^2 R(U) \left( \nabla U , \nabla I_h^p(U^{q-1})\right)_{0,T} 
+ C h_T^{\beta} q^{p+1}. 
\label{gesamt_int_error}
\end{eqnarray}
Inserting this into $(\ref{weak_problem4})$ we obtain
\begin{equation}
\begin{split}
& \int_\Omega \left( U_-^N \right)^q \, d\sx  
-C q^{p+2} \sum_{n=0}^{N-1} \sum_{T \in \T_h} 
h_T^2  R(U) \left( \nabla U , \nabla I_h^p(U^{q-1})\right)_{0,T} \\
&+q \sum_{n=0}^{N-1} \sum_{T \in \T_h^n}  \hat{\eps}(U) \left(\nabla U, \nabla I_h^p(U^{q-1}) \right)_{0,T}
+ \left(\delta(U) L(U), \stf'(U) \nabla (U^{q-1})\right)_{0,T}\\
&+\sum_{i=4}^5 E_i(f,U^q/q,U,1)
\leq \int_\Omega \left( u_0 \right)^q d\sx
+q F(f,U^q/q,U,1) +C h_T^{\beta} q^{p+2}.
\end{split}
\label{weak_problem5}
\end{equation}
To proceed with the treatment of interpolation error, it is necessary to require
that $C q^{p+2} \leq h^{-\beta},$ where $0 < \beta < 1/2.$
The upper limit of $\beta$ is introduced due to convergence reasons,
cf. \cite[Lemma 3.2]{Jaffre.Johnson.ea_Convergenceofdiscontinuous_1995}.
However, this restriction on $q$ does not prevent us
to finish this proof by letting $q \to \infty.$

Moreover,  we have
\begin{equation}
\begin{split}
& \int_\Omega \left( U_-^N \right)^q  d\sx  
+q \sum_{n=0}^{N-1} \sum_{T \in \T_h^n} 
 \left( \delta(U)L(U), \stf'(U) \nabla (U^{q-1})\right)_{0,T} \\
&+\sum_{i=4}^5 E_i(f,U^q/q,U,1)\leq \int_\Omega \left( u_0 \right)^q d\sx
+q F(f,U^q/q,U,1) +C.
\end{split}
\label{weak_problem6}
\end{equation}
Note that
\begin{eqnarray*}
& & E_4(f,U^q/q,U,1)+E_5(f,U^q/q,U,1)\\
& \geq& \sum_{n=0}^{N-1} \sum_{\tau \in \Lambda_{n,n+1}}
\int_\tau \left( C_0^{\partial \Omega} -\frac{1}{2}
\int_0^1 \left|\stf'(g_D +\xi  \llbracket U \rrbracket ) \cdot n^+\right|\,d\xi  \right)
\llbracket U^{q-1}\rrbracket  \llbracket U\rrbracket  \, ds \notag, 
\end{eqnarray*}
and thus, by Young's inequality,
\begin{eqnarray}
& & F(f,U^q/q,U,1) \notag \\
&\leq& \sum_{n=0}^{N-1} \sum_{\tau \in \Lambda_{n,n+1}}
\int_\tau \left( \frac{1}{2} \int_0^1 | \stf'(g_D + \xi \llbracket U \rrbracket) \cdot n^+| \, d\xi
+C_0^{\partial \Omega} \right)\, | \llbracket U \rrbracket g_D^{q-1}| \, \, ds \notag \\
 &\underset{(\ref{r01})}{\leq} &\sum_{n=0}^{N-1}
\sum_{\tau \in \Lambda_{n,n+1}} \int_\tau \frac{1}{2} C_0^{\partial \Omega}
\, \llbracket U^{q-1} \rrbracket \llbracket U \rrbracket \, ds 
+ \frac{1}{2} C_{0}^{\partial \Omega}
(q-1)q^{-\frac{q}{q-1}} 3^{\frac{q}{q-1}} 2^{\frac{q-2}{q-1}} g_D^q  \, ds \notag \\
&=&  F_1(f,U^q/q,U,1) +F_2(f,g_D^q/q,g_D,1). \label{f45_uu} 
\end{eqnarray}
So we obtain
\begin{equation*}
E_4(f,U^q/q,U,1)+E_5(f,U^q/q,U,1)-F_1(f,U^q/q,U,1) \geq 0.
\end{equation*}
Altogether we get that
\begin{equation}
\begin{split}
 \int_\Omega \left( U_-^N \right)^q  d\sx  
&+q \sum_{n=0}^{N-1} \sum_{T \in \T_h^n} 
 \left( \delta(U)L(U), \stf'(U) \nabla (U^{q-1})\right)_{0,T} \\
&\leq \int_\Omega \left( u_0 \right)^q d\sx
+ C_0 (q-1) \left( 2q \right)^{-\frac{1}{q-1}}3^{\frac{q}{q-1}} \|g_D \|_{0,q,\Sigma_T}^q+ C.\\
\end{split}
\label{weak_problem7}
\end{equation}
By repeating the arguments given at the end of the previous section,
we summarize that
\begin{equation}
\underset{t \geq 0}{\sup} \|U(t,\cdot)\|_{0,q,\Omega} \leq C^{\frac{1}{q}} \left(
\|u_0\|_{0,q,\Omega}+ 3^{\frac{1}{q-1}}\left(C_0 q \right)^{\frac{1}{q}}    \|g_D \|_{0,q,\Sigma_T}
+ C^{\frac{1}{q}}\right)
\label{L_q_abschaetzung}
\end{equation}
for $4 \leq q \leq Ch^{-\frac{\beta}{p+2}}.$
Finally, using an inverse inequality we have that
\begin{equation*}
\begin{split}
\|U\|_{0,\infty,Q_T} &\leq \left( Cq h^{-1} \right)^{\frac{\sd}{q}} \underset{t \geq 0}{\sup} \|U(t,\cdot)\|_{0,q,\Omega} 
\leq \left( C h^{-1-\frac{\beta}{p+2}}\right)^{\frac{\sd}{q}} \underset{t \geq 0}{\sup} \|U(t,\cdot)\|_{0,q,\Omega} \\
&= C^{\frac{\sd}{q}} \, \exp\left(C \sd q^{-1} \ln \left(h^{-1}\right)\right) \underset{t \geq 0}{\sup} \|U(t,\cdot)\|_{0,q,\Omega}.
\end{split}
\end{equation*}
Setting $q= C h^{-\frac{\beta}{p+2}}$ we get 
\begin{equation}
\begin{split}
\|U\|_{0,\infty,Q_T} &\leq C^{\sd h^{\beta/(p+2)}} \, \exp\left(C \sd h^{\beta/(p+2)}
\ln \left(h^{-1}\right)\right) \underset{t \geq 0}{\sup} \|U(t,\cdot)\|_{0,q,\Omega} \\
& \leq  C \underset{t \geq 0}{\sup} \|U(t,\cdot)\|_{0,q,\Omega}
\end{split}
\label{u_infty_q}
\end{equation}
in the case of $h\leq1,$ which concludes the proof.
\end{proof}

\begin{corollary}
Under the assumptions of Theorem \ref{Th:Linfty_main}, the estimate
\begin{equation}
\|U\|_{0,\infty,Q_T} \leq  
 \|u_0\|_{0,\infty,\Omega} +  \|g_D \|_{0,\infty,\Sigma_T} +1 
\label{u_infty_h_to_0}
\end{equation}
holds for $h \to 0.$
\end{corollary}
\begin{proof}
Using $(\ref{L_q_abschaetzung})$ and $(\ref{u_infty_q})$ with
$h=C q^{-(p+2)/\beta} \to 0,$ the statement immediately follows.
\end{proof}

\begin{remark}
The $L^\infty(L^\infty)$-boundedness of DG($p$)-solutions with $p>1$
was also considered in \cite{Szepessy_Convergenceofstreamline_1991}.
However, since an inequality of the form $(\ref{nabla_qm1_ungl})$
was proved only for the case $p=1,$ the shock-capturing term was realized
on finer auxiliary triangulations using polynomials of first degree.
Hence, the bound $C^{q}q^{p+2} \leq h^{-\beta}$ is necessary,
which is true for $q \leq C \ln(1/h),$
cf. \cite[(3.16)]{Szepessy_Convergenceofstreamline_1991}.
This gives
\begin{equation*}
\|U\|_{0,\infty,Q_T} \leq 
 C \underset{t \geq 0}{\sup} \|U(t,\cdot)\|_{0,q,\Omega},
\end{equation*}
for $h \to 0,$ where $C \neq 1.$
\label{bem:L_infty}
\end{remark}

\begin{corollary}
Under the assumptions of Theorem \ref{Th:Linfty_main}, the estimate
\begin{equation}
\|U\|_{0,\infty,Q_T} \leq  
 \|u_0\|_{0,\infty,\Omega} + \|g_D \|_{0,\infty,\Sigma_T} \quad \forall h>0
\label{u_infty_p_0}
\end{equation}
holds for $p=0.$
\end{corollary}
\begin{proof}
Let $p=0.$ Then we get $\left( b(U,U^{q-1})-b(U,I_h^p(U^{q-1})) \right)=0.$
Consequently, there is no need for the bound $q \leq Ch^{-\frac{\beta}{p+2}},$
and we can conclude with letting $q \to \infty$ in $(\ref{L_q_abschaetzung})$
and $(\ref{u_infty_q}).$
\end{proof}

\section{Proof of Lemma \ref{nabla_qm1}}
Until now there is no proof of an inequality like $(\ref{nabla_qm1_ungl})$
for $q \neq 2.$ To the best of our knowledge only special cases for
linear ansatz functions on triangles respectively tetrahedrons are available,
cf. \cite[Lemma 4.2]{Szepessy_Convergenceofshock-capturing_1989} and
\cite[Lemma 3.3]{Szepessy_Convergenceofstreamline_1991}.
Moreover, the constant in these references depends on $q.$   

Using the theory of numerical ranges for bounded linear operators in Banach spaces,
we  are able to prove this inequality under rather weak assumptions.
More precisely,
the local stiffness matrix of the shock-capturing term has to be
symmetric positively definite and an eigenvector $(1,\cdots,1)^T$
with an unique eigenvalue zero.

First of all we need some further notation and definitions about
the numerical range.
\begin{definition}
Let $(X,\|\cdot\|)$ be a normed vector space, let $S(X)$ be the unit sphere and
denote by $X'$ the dual space of $X.$ 
For each bounded linear operator $A$ on $X,$
\begin{equation}
W(A,\, \|\cdot\|)=\left\{ f(Ax):\, (x,f) \in \Pi \right\}
\label{def_wertebereich}
\end{equation}
with
$\Pi=\left\{ (x,f) \in S(X) \times S(X'): \, f(x)=1 \right\}$
is called the spatial numerical range.
\end{definition}
\pagebreak
\begin{remark}{\ } %
\begin{enumerate}
\item Note that the following definition is equivalent to (\ref{def_wertebereich}):
\begin{equation}
W(A,\, \|\cdot\|)=\left\{ f(Ax) :\, f(x)=\|x\| \|f\|=1 \right\},
\label{def_wertebereich_bauer}
\end{equation}
since 
\begin{equation*}
\|f\|=\underset{\|x\|=1}{\sup} |f(x)| \leq \|f\| \underbrace{\|x\|}_{=1}=1.
\end{equation*}
\item In contrast to the spectrum $\sigma(A),$ the spatial numerical range
$W(A,\, \|\cdot\|)$ depends on the norm $\|\cdot \|.$ 
\item Let  $\l^{q_1}(n)=(\R^n,\|\cdot \|_{\l^{q_1}})$ be the normed vector space.
Due to $f(x)=\sum_{i=1}^{n} x_i f(e_i)=x^T y_f$ and $\l^{q_2}(n)$ for $1/{q_1}+1/{q_2}=1,$
which is norm-isomorphic to $\l^{q_1}(n)',$ the identity $\|f\|=\|y_f\|_{\l^{q_2}}$ is valid. 
Therefore we get 
\begin{equation}
W(A,\, \|\cdot\|_{\l^{q_1}})=\left\{ x^T A y_f :\, x^T y_f=\|x\|_{\l^{q_1}} \|y_f\|_{\l^{q_2}}=1 \right\},
\label{def_wertebereich_bauer2}
\end{equation}
which is the identity case of Hölder's inequality, cf. \cite{Bauer_fieldofvalues_1962}.
In the case ${q_1}=2,$ $W(A,\, \|\cdot\|_{\l^2})=W(A)$ is the numerical range
in a Hilbert space due to Toeplitz \cite{Toeplitz_DasalgebraischeAnalogon_1918}. 
\item Unlike to $W(A),$ the spatial numerical range is not necessary convex,
cf. \cite[S. 357]{Nirschl.Schneider_Bauerfieldsof_1964}.
\end{enumerate}
\end{remark}

If we interpret the Matrix $A$ as an element of a normed algebra $\mathcal{A}$
with an identity element, we can define a second numerical range.
For further details we refer to
\cite[S. 15]{Bonsall.Duncan_Numericalrangesof_1971}.

\begin{definition}
Let $\mathcal{A}$ be a normed algebra,
$S(\mathcal{A})=\left\{ x \in \mathcal{A} : \|x\|=1 \right\}$ the unit sphere
and $\mathcal{A}'$ the dual space of $\mathcal{A}.$
For  $x \in \mathcal{A}$ let
\begin{equation}
D(\mathcal{A},x)=\left\{ f \in \mathcal{A}':f(x)=1=\|f\| \right\}.
\label{dual_algebra_element} 
\end{equation}
We define 
the algebraic numerical range by
\begin{equation}
V_{\mathcal{A}}(a,\|\cdot\|)
=\cup \left\{ V_{\mathcal{A}}(a,x,\|\cdot\|) : x \in S(\mathcal{A}) \right\},
\label{num_bereich_alg2}
\end{equation}
where
\begin{equation}
V_{\mathcal{A}}(a,x,\|\cdot\|)=\left\{ f(ax): f \in D(\mathcal{A},x) \right\}.
\label{num_bereich_alg1} 
\end{equation} 
\end{definition}

Notice that for the algebraic numerical range it is sufficient to consider only
the identity element.
\begin{lemma}
\begin{equation*}
V_{\mathcal{A}}(a,\|\cdot\|)=V_{\mathcal{A}}(a,1,\|\cdot\|),\, a \in \mathcal{A}.
\end{equation*}
\end{lemma}
\begin{proof}
\cite[Lemma 2.2]{Bonsall.Duncan_Numericalrangesof_1971}.
\end{proof}
Next, we recall two well-known results about numerical ranges and the connection
to the spectrum $\sigma(A).$
\begin{lemma}
\begin{equation}
\conv W(A,\|\cdot\|) =V_{\mathcal{A}}(a,\|\cdot \|).
\label{id_W_V}
\end{equation}
\end{lemma}
\begin{proof}
\cite[S. 84]{Bonsall.Duncan_Numericalrangesof_1971}.
\end{proof}
\begin{theorem}[Vidav]
Let $a \in \mathcal{A}$ be a Hermitian element, i.e.,
$V_{\mathcal{A}}(a,\|\cdot\|) \subset \R.$ Then we have
\begin{equation}
\conv \sigma(a)=V_{\mathcal{A}}(a,\|\cdot\|).
\label{vidav}
\end{equation}
\end{theorem}
\begin{proof}
\cite[Corollary 5.11]{Bonsall.Duncan_Numericalrangesof_1971}.
\end{proof}
Further we have a corollary which will help us to prove Lemma \ref{nabla_qm1}.
\begin{corollary}
Let $A$ be a symmetric, positively semidefinite matrix. Then we have
\begin{equation}
W(A,\|\cdot\|) \subseteq [0,\lambda_{\max}(A)],
\label{folg_WA}
\end{equation}
where $\lambda_{\max}(A)$ denotes the largest eigenvalue of $A.$
\end{corollary}

\begin{proof}[Proof of Lemma \ref{nabla_qm1}]
First, consider  $(\ref{nabla_qm1_ungl})$ on the reference element. 
Obviously, the inequality is valid for $v=\const.$ Let $v \neq \const$ be given. 
Consider a decomposition of 
$\mathbb{Q}_p(\hat{T})=V^0(\hat{T}) \oplus V(\hat{T})$ with 
\begin{eqnarray*}
V^0(\hat{T})&=&\{v \in \mathbb{Q}_p(\hat{T}):v=\const\}\\
      &=&\{v \in \mathbb{Q}_p(\hat{T}):\int_{\hat{T}} \nabla w \cdot \nabla v\, dx =0,
\, w \in \mathbb{Q}_p(\hat{T})\}.
\end{eqnarray*}
Let $V_{\mathcal{N}}$ and $V_{N}^0,$ $V^0_{\mathcal{N}} \oplus V_{\mathcal{N}}=\R^{\ndofp}$
be the coefficient spaces.
Therefore, due to the definition of Lagrange finite elements we have that
\begin{equation}
V_{\mathcal{N}}^0=\myspan\{(1,\dots,1)^T\}, \, \dim V_{\mathcal{N}}^0=1. 
\label{dimVN0}
\end{equation}
Let $\mathcal{N}$ denote the Lagrange nodes and $\ndofp$ the number of degrees of freedom.
Moreover, defining $\nabla \varphi = (\nabla \varphi_{1},\dots,\nabla \varphi_{\ndofp})^T$,
$v_{\mathcal{N}}=(v(x))_{x \in \mathcal{N}},$
$v_{\mathcal{N}}^{q-1}\linebreak
=(v^{q-1}(x))_{x \in \mathcal{N}},$
we have that
\begin{equation*}
\frac{
\left( \nabla v, \nabla I_h^p(v^{q-1})\right)_{0,\hat{T}}
}{\|v_{\mathcal{N}}\|_{\l^q}^q
}
=\frac{v_{\mathcal{N}}^T \hat{A} v_{\mathcal{N}}^{q-1}}{v_{\mathcal{N}}^T v_{\mathcal{N}}^{q-1}}
= v_{\mathcal{N}}^T \hat{A} v_{\mathcal{N}}^{q-1}
\end{equation*}
if $v_{\mathcal{N}}^T v_{\mathcal{N}}^{q-1}=1.$ 
Notice that, due to the homogeneity of the quotient, such a norming is always possible. 
Using the fact that $\hat{A}$ is a symmetric, positively semidefinite matrix and 
$1=v_{\mathcal{N}}^T v_{\mathcal{N}}^{q-1}=\|v_{\mathcal{N}}\|_{\l^q} \|v_{\mathcal{N}}^{q-1}\|_{\l^{q/(q-1)}},$
we obtain 
\begin{equation*}
v_{\mathcal{N}}^T \hat{A} v_{\mathcal{N}}^{q-1} \in W(\hat{A},\|\cdot\|_{\l^q})
\underset{(\ref{folg_WA})}{\subseteq} [0,\lambda_{\max}(\hat{A})]. 
\end{equation*}
Now, it is natural to ask whether $v_{\mathcal{N}}^T \hat{A} v_{\mathcal{N}}^{q-1}$
is bounded from zero independent of $q.$ 
To see this, let us suppose that the eigenvalues are ordered in increasing manner
\begin{equation*}
0=\lambda_1 < \lambda_2 \leq \dots \leq \lambda_{\ndofp}.
\end{equation*}
Since $V_{\mathcal{N}}^0$ is the eigenspace of $\lambda_1,$ we have
\begin{equation}
V^0_{\mathcal{N}} \perp V_{\mathcal{N}},
\label{eq:Qp_perp}
\end{equation}
and we can write $\hat{A}$ in terms of a sum of dyadic products
of eigenvectors $\xi_i, 1\leq i\leq \ndofp$ 
\begin{equation*}
v_{\mathcal{N}}^T \hat{A} v_{\mathcal{N}}^{q-1}  
=v_{\mathcal{N}}^T \sum_{i=1}^{\ndofp} \lambda_i \xi_i \xi_i^T v_{\mathcal{N}}^{q-1}  
=\sum_{i=1}^{\ndofp} \lambda_i v_{\mathcal{N}}^T  \xi_i \xi_i^T v_{\mathcal{N}}^{q-1}  
=\sum_{i=2}^{\ndofp} \lambda_i v_{\mathcal{N}}^T  \xi_i \xi_i^T v_{\mathcal{N}}^{q-1}. 
\end{equation*}
Again, the inclusion $(\ref{folg_WA})$ yields
$v_{\mathcal{N}}^T  \xi_i \xi_i^T v_{\mathcal{N}}^{q-1} \in [0,1].$
On the other hand, supposing
$v_{\mathcal{N}} \neq \const,$
we obtain
\begin{equation}
v_{\mathcal{N}}, v_{\mathcal{N}}^{q-1} \perp \myspan\{\xi_1\}=V_{\mathcal{N}}^0.
\label{eq:v_vp_perp_xi1}
\end{equation}
Then, the boundedness of $v_{\mathcal{N}}^T \hat{A} v_{\mathcal{N}}^{q-1},$ i.e.
\begin{equation*}
v_{\mathcal{N}}^T \hat{A} v_{\mathcal{N}}^{q-1}  
\geq 
\lambda_2 \sum_{i=1}^{\ndofp} v_{\mathcal{N}}^T  \xi_i \xi_i^T v_{\mathcal{N}}^{q-1} 
- \lambda_2 \underbrace{v_{\mathcal{N}}^T  \xi_1}_{=0} \xi_1^T v_{\mathcal{N}}^{q-1} 
\geq \lambda_2 v_{\mathcal{N}}^T v_{\mathcal{N}}^{q-1},
\end{equation*}
implies that
\begin{equation}
\frac{
\left( \nabla v, \nabla I_h^p(v^{q-1})\right)_{0,\hat{T}}
}{\|v_{\mathcal{N}}\|_{\l^q}^q
}
\geq \lambda_2,\quad v\neq \const.
\label{koerziv1}
\end{equation}
Now, standard estimates give the proof for the reference element
\begin{eqnarray*}
\|v\|_{0,\infty,\hat{T}}^{q-2} \|\nabla v\|_{0,2,\hat{T}}^2 
&\leq& \lambda_{\max}(\hat{A}) \|v^{q-2}\|_{0,\infty,\hat{T}} \|v_{\mathcal{N}}\|_{\l^2}^2 \\
&\leq& \lambda_{\max}(\hat{A})\Lambda_p \|v_{\mathcal{N}}\|_{\l^\infty}^{q-2} \|v_{\mathcal{N}}\|_{\l^2}^2 \\
&\leq& \lambda_{\max}(\hat{A})\Lambda_p \|v_{\mathcal{N}}\|_{\l^q}^{q-2} \|v_{\mathcal{N}}\|_{\l^2}^2 \\
&\leq& \left( \ndofp \right)^{1-2/q} \lambda_{\max}(\hat{A})\Lambda_p \|v_{\mathcal{N}}\|_{\l^q}^{q} \\
&\leq& \ndofp \frac{\lambda_{\max}(\hat{A})}{\lambda_2}\Lambda_p 
\left( \nabla v, \nabla I_h^p(v^{q-1})\right)_{0,\hat{T}}, 
\end{eqnarray*}
where $\Lambda_p=\|\sum_{i=1}^{\ndofp}|\varphi_i| \|_{0,\infty,\hat{T}}$
is the Lebesgue constant.

Finally, we want to show the result for an affine decomposition
$\left\{ \mathcal{T}_h^n \right\}_{h>0}$ with 
\begin{equation*}
F_T:\hat{T} \ni \hat{x} \mapsto J_T \hat{x} +b_T =x \in T
\qquad \forall T \in \left\{ T_h^n \right\}_{h>0}
\end{equation*}
and therefore
\begin{equation}
\nabla u(x) =J_T^{-T} \hat{\nabla} \hat{u}(\hat{x}), \quad \hat{u}=u \circ F_T.
\label{aff_grad}
\end{equation}
Due to the spectral decomposition of 
\begin{equation*}
K=\left( J_T^T J_T \right)^{-1}=\sum_{l=1}^d \mu_l \psi_l \psi_l^T,
\end{equation*}
we have that
\begin{eqnarray*}
\int_T \nabla v \cdot \nabla I_h^p(v^{q-1})\, dx
&=& v_{\mathcal{N}}^T \left( \int_T \nabla \varphi_j \cdot \nabla \varphi_i \, dx \right)_{i,j} v_{\mathcal{N}}^{q-1}\\
&=&v_{\mathcal{N}}^T \left( \int_{\hat{T}} (\hat{\nabla} \hat{\varphi}_j)^T
K \hat{\nabla} \hat{\varphi}_i | \det(J_T)|  \, d\hat{x} \right)_{i,j} v_{\mathcal{N}}^{q-1}\\
&=&| \det(J_T)| \sum_{l=1}^d \mu_i  v_{\mathcal{N}}^T
\left( \int_{\hat{T}} (\hat{\nabla} \hat{\varphi}_j)^T \psi_l \psi_l^T
\hat{\nabla} \hat{\varphi}_i   \, d\hat{x} \right)_{i,j}
v_{\mathcal{N}}^{q-1}\\
&\underset{(\ref{folg_WA})}{\geq} & | \det(J_T)|  \mu_{\min}  v_{\mathcal{N}}^T \hat{A}v_{\mathcal{N}}^{q-1}\\
&=& | \det(J_T)|  \|J_T\|_{\l^2}^{-2}  v_{\mathcal{N}}^T \hat{A}v_{\mathcal{N}}^{q-1}
\end{eqnarray*}
and
\begin{equation*}
|\hat{\nabla} \hat{v}|_{0,2,\hat{T}}^2 \geq |\det(J_T)|^{-1} \|J_T^{-1}\|_{\l^2}^{-2} |\nabla v|_{0,2,T}^2.
\end{equation*}
Now, the shape regularity property 
\begin{equation*}
\|J_T\|_{\l^2} \|J_T^{-1}\|_{\l^2} \leq C
\end{equation*}
completes the proof.
\end{proof}

\section{Summary}
In this paper, we considered a DG-method based on polynomials of degree $p\geq 0$
for hyperbolic scalar conservation laws.
This method was introduced for the pure Cauchy problem in
\cite{Jaffre.Johnson.ea_Convergenceofdiscontinuous_1995}.
We extended the formulation for hyperbolic conservation laws
with initial and boundary conditions.
Moreover, we presented a proof of the uniform boundedness
of the discrete solution in the $L^\infty(L^\infty)$-norm.
The analysis is based on arguments demonstrated
in \cite{Szepessy_Convergenceofstreamline_1991} which are valid for $p=1.$
It turned out that the use of numerical ranges for bounded linear operators
in Banach spaces allows to generalize this result to the case $p>1.$
Future work will be devoted to the convergence of the DG-method
for the initial-boundary value problem.

\addcontentsline{toc}{chapter}{Literaturverzeichnis}


\begin{thebibliography}{MNRR96}

\bibitem[Bau62]{Bauer_fieldofvalues_1962}
F.~L. Bauer.
\newblock On the field of values subordinate to a norm.
\newblock {\em Numerische Mathematik}, 4:103--113, 1962.

\bibitem[BD71]{Bonsall.Duncan_Numericalrangesof_1971}
F.~F. Bonsall and J.~Duncan.
\newblock {\em Numerical ranges of operators on normed spaces and of elements
  of normed algebras}, volume~2 of {\em London Mathematical Society Lecture
  Note Series}.
\newblock Cambridge University Press, London, 1971.

\bibitem[BlRN79]{Bardos.Roux.ea_Firstorderquasilinear_1979}
C.~Bardos, A.~Y. le~Roux, and J.-C. N{\'e}d{\'e}lec.
\newblock First order quasilinear equations with boundary conditions.
\newblock {\em Communications in Partial Differential Equations},
  4(9):1017--1034, 1979.

\bibitem[EG04]{Ern.Guermond_Theoryandpractice_2004}
A.~Ern and J.~Guermond.
\newblock {\em Theory and practice of finite elements}.
\newblock Springer-Verlag, New York, 2004.

\bibitem[JJS95]{Jaffre.Johnson.ea_Convergenceofdiscontinuous_1995}
J.~Jaffr{\'e}, C.~Johnson, and A.~Szepessy.
\newblock Convergence of the discontinuous galerkin finite element method for
  hyperbolic conservation laws.
\newblock {\em Math. Models Methods Appl. Sci.}, 5(3):367--386, 1995.

\bibitem[JSH90]{Johnson.Szepessy.ea_convergenceofshock-capturing_1990}
C.~Johnson, A.~Szepessy, and P.~Hansbo.
\newblock On the convergence of shock-capturing streamline diffusion finite
  element methods for hyperbolic conservation laws.
\newblock {\em Math. Comp.}, 54(189):107--129, 1990.

\bibitem[MNRR96]{Malek.Nevcas.ea_WeakandMeasure-valued_1996}
J.~Malek, J.~Ne\v{c}as, M.~Rokyta, and M.~R\r{u}\v{z}i\v{c}ka.
\newblock {\em Weak and Measure-valued Solutions to Evolutionary PDEs},
  volume~13 of {\em Applied Mathematics and Mathematical Computation}.
\newblock Chapman \& Hall, 1996.

\bibitem[NS64]{Nirschl.Schneider_Bauerfieldsof_1964}
N.~Nirschl and H.~Schneider.
\newblock The {B}auer fields of values of a matrix.
\newblock {\em Numer. Math.}, 6:355--365, 1964.

\bibitem[Sze89a]{Szepessy_Convergenceofshock-capturing_1989}
A.~Szepessy.
\newblock Convergence of a shock-capturing streamline diffusion finite element
  method for a scalar conservation law in two space dimensions.
\newblock {\em Math. Comp.}, 53(188):527--545, 1989.

\bibitem[Sze89b]{Szepessy_Measure-valuedsolutionsof_1989}
A.~Szepessy.
\newblock Measure-valued solutions of scalar conservation laws with boundary
  conditions.
\newblock {\em Arch. Rational Mech. Anal.}, 107(2):181--193, 1989.

\bibitem[Sze91]{Szepessy_Convergenceofstreamline_1991}
A.~Szepessy.
\newblock Convergence of a streamline diffusion finite element method for
  scalar conservation laws with boundary conditions.
\newblock {\em RAIRO Mod\'el. Math. Anal. Num\'er.}, 25(6):749--782, 1991.

\bibitem[Toe18]{Toeplitz_DasalgebraischeAnalogon_1918}
O.~Toeplitz.
\newblock Das algebraische {A}nalogon zu einem {S}atze von {F}ej\'er.
\newblock {\em Math. Z.}, 2(1-2):187--197, 1918.

\end{thebibliography}
\end{document}